\documentclass[10pt]{article}

\usepackage{amsmath}
\makeatletter
\newfont{\footsc}{cmcsc10 at 8truept}
\newfont{\footbf}{cmbx10 at 8truept}
\newfont{\footrm}{cmr10 at 10truept}
\renewcommand{\ps@plain}{}
\makeatother 

\title{\bf The Smallest Degree Sum that Yields Potentially $C_k$-graphical
Sequences}

\author{ Chunhui Lai \thanks{ Project Supported by NNSF of China(10271105), NSF of Fujian,
Science and Technology Project of Fujian, Fujian Provincial
Training Foundation for "Bai-Quan-Wan Talents Engineering" ,
Project of Fujian Education Department and Project of Zhangzhou
Teachers College.}\\
\small Department of Mathematics, Zhangzhou Teachers College,\\
\small Zhangzhou, Fujian 363000, P. R. of CHINA\\
\small \texttt{ zjlaichu@public.zzptt.fj.cn }}
\date{ Journal of Combinatorial Mathematics and Combinatorial Computing  49 (2004), 57-64.}
\begin{document}
\maketitle

\begin{abstract}
In this paper we consider a variation of the classical
Tur\'{a}n-type extremal problems. Let $S$ be an $n$-term graphical
sequence, and $\sigma(S)$ be the sum of the terms in $S$. Let $H$
be a graph. The problem is to determine
 the smallest even $l$ such that any $n$-term graphical sequence $S$ having
$\sigma(S)\ge l$ has a realization containing $H$ as a subgraph.
Denote this value $l$ by $\sigma(H,\ n)$.\ We show
$\sigma(C_{2m+1},\ n)=m(2n-m-1)+2$,\ for $m\ge 3$,\ $n\ge 3m$;\
$\sigma(C_{2m+2},\ n)=m(2n-m-1)+4$,\ for $m\ge 3,\ n\ge 5m-2$.
 \end{abstract}

  \par

\section{Introduction}
\par

  If $S=(d_1,d_2,...,d_n)$ is a sequence of
non-negative integers, then it is called  graphical if there is a
simple graph $G$ of order $n$, whose degree sequence ($d(v_1 ),$
$d(v_2 ),$ $...,$ $d(v_n )$) is precisely $S$. If $G$ is such a
graph then $G$ is said to realize $S$ or be a realization of $S$.
A graphical sequence $S$ is potentially $H$ graphical if there is
a realization of $S$ containing $H$ as a subgraph, while $S$ is
forcibly $H$ graphical if every realization of $S$ contains $H$ as
a subgraph. Let $\sigma(S)=d(v_1 )+d(v_2 )+... +d(v_n ),$ and
$[x]$ denote the largest integer less than or equal to $x$. If $G$
and $G_1$ are graphs, then $G\cup G_1$ is the disjoint union of
$G$ and $G_1$. If $G = G_1$, we abbreviate $G\cup G_1$ as $2G$.
Let $K_k$, and $C_k$ denote a complete graph on $k$ vertices, and
a cycle on $k$ vertices, respectively.\par

Given a graph $H$, what is the maximum number of edges of a graph
with $n$ vertices not containing $H$ as a subgraph? This number is
denoted $ex(n,H)$, and is known as the Tur\'{a}n number. This
problem was proposed for $H = C_4$ by Erd\"os [2] in 1938 and in
general by Tur\'{a}n [11]. In terms of graphic sequences, the
number $2ex(n,H)+2$ is the minimum even integer $l$ such that
every $n$-term graphical sequence $S$ with $\sigma (S)\geq l $ is
forcibly $H$ graphical. Here we consider the following variant:
determine the minimum even integer $l$ such that every $n$-term
graphical sequence $S$ with $\sigma(S)\ge l$ is potentially $H$
graphical. We denote this minimum $l$ by $\sigma(H, n)$. Erd\"os,\
Jacobson and Lehel [3] showed that $\sigma(K_k, n)\ge
(k-2)(2n-k+1)+2$ and conjectured that equality holds. They proved
that if $S$ does not contain zero terms, this conjecture is true
for $k=3,\ n\ge 6$. The conjecture is confirmed in [4],[7],[8],[9]
and [10].
 \par
 Gould,\ Jacobson and
Lehel [4] also proved that  $\sigma(pK_2, n)=(p-1)(2n-2)+2$ for
$p\ge 2$; $\sigma(C_4, n)=2[{{3n-1}\over 2}]$ for $n\ge 4$. Lai
[5, 6] proved that $\sigma(C_5, n)=4n-4$ for $n\ge 5$, and
$\sigma(C_6, n)=4n-2$ for $n\ge 7$, $\sigma(C_{2m+1}, n)\ge
m(2n-m-1)+2$,\ for $n\ge 2m+1, m\ge 2$, $\sigma (C_{2m+2}, n) \ge
m(2n-m-1)+4$, for $n\ge 2m+2, m\ge 2$, $\sigma (K_4-e,
n)=2[{{3n-1}\over 2}]$ for $n\ge 7$.\ In this paper we prove that
$\sigma (C_{2m+1}, n)=m(2n-m-1)+2$,\
 for $n\ge 3m$,\ $m\ge 3$;\ $\sigma(C_{2m+2}, n)=m(2n-m-1)+4$ for $n\ge 5m-2,\ m\ge 3$.\par
 \section{ Main results.} \par
 {\bf Theorem\
1.} Let $k\ge 4$. Let $S$ be a potentially $C_k$-graphical n-term
sequence. If there exists $x\notin C_k, w\in C_k$ such that
$d(x)\ge [{k\over 2}]+1$,\ $d(w)\ge 3$. Then $S$ has a realization
containing a $C_{k+1}$.\par Assume $C_k$ is $w_1 w_2\cdots w_k
w_1$. Let $w_{k+i}=w_i$.\ We first give the following three
results.
\par {\bf Lemma (a)} For any $x\notin C_k$,\ if there
is $w_r, w_{r+1}$ such that $w_r x, w_{r+1}x\in E(G)$, then $G$
contains a $C_{k+1}:\ w_1 w_2\cdots w_r x w_{r+1}\cdots w_k
w_1$.\par {\bf Lemma(b)} For any $ x, y\notin C_k, xy\in E(G)$, if
there is $w_r$ such that $w_r x\in E(G)$,\ $w_r y\notin E(G)$,
then $S$ has a realization containing a $C_{k+1}$. (We see the
edge $w_{r+1} x$ is not in $G$ or a $C_{k+1}$ would
 exist, but then the edge interchange which removes the  edges
$w_r w_{r+1}$ and $xy$ and inserts the  edges $w_{r+1} x$ and $w_r
y$ produces a realization containing a $C_{k+1}:\ w_1 w_2\cdots
w_r x w_{r+1}\cdots w_k w_1$)\par {\bf Lemma(c)} For any $ x,
y\notin C_k$, $xy\in E(G)$, if there is $w_r, w_{r+2}$ such that
$w_r x,\ w_{r+2}x\in E(G)$, then $S$ has a realization containing
a $C_{k+1}$. (If $w_{r+2} y\notin E(G)$, then by Lemma(b), $S$ has
a realization containing a $C_{k+1}$. Otherwise, $w_{r+2} y\in
E(G)$ and so $G$ contains a $C_{k+1}: w_1 w_2\cdots w_r xy$
$w_{r+2} w_{r+3}\cdots w_k w_1$)
\par  {\bf Proof of theorem\ 1.} Assume every
realization of $S$ does not contain a $C_{k+1}$. By Lemma(a), $x$
is  adjacent to at most $[{k\over 2}]$ vertices of $C_k$. Since
$d(x)\ge [{k\over 2}]+1$ there exists $x_1\notin C_k$ such that
$xx_1\in E(G)$. Thus, by  Lemma(c), $x$ is adjacent to at most
$[{k\over 3}]$ vertices of $C_k$. Note that  $[{k\over 3}] \le
[{k\over 2}]-1$ since $k\ge 4$. Hence there is $x_2\notin C_k$,
$x_2\ne x_1$, such that $x x_2\in E(G)$.
\par Case 1.\ Suppose that there is
$w_i\in C_k$ such that $w_i x\in E(G)$. By Lemma(b), $w_i x_1, w_i
x_2\in E(G)$. By Lemma(a), $w_{i+1} x, w_{i+1} x_1, w_{i+1}
x_2\notin E(G)$. By Lemma(c) $w_{i+2} x, w_{i+2} x_1, w_{i+2}
x_2\notin E(G)$. Then the edge interchange which removes the
 edges $w_{i+1} w_{i+2}$ and $x x_2$ and inserts the
 edges $w_{i+2} x$ and $w_{i+1} x_2$ produces a realization
containing a $C_{k+1}:\ w_1 w_2$ $\cdots w_i x_1 x w_{i+2}
w_{i+3}$ $\cdots w_k w_1$. This is a contradiction.
\par Case 2.\ Suppose
for any $ w_i\in C_k, w_i x\notin E(G)$.\ Since $d(x)\ge [{k\over
2}]+1\ge 2+1=3$, hence there is $x_3\notin C_k$, $x_3\ne x_1$,
$x_3\ne x_2$ such that $x x_3\in E(G)$. By Lemma(b), $w_i x_1,w_i
x_2,w_i x_3\notin E(G)$. Since there is $w\in C_k$ such that
$d(w)\ge 3$, then there is $x_4$
 such that $w x_4\notin E(C_k)$, $wx_4\in E(G)$. By Lemma(b), $x_4$
 is not one of $ x_1$,$ x_2$,$ x_3$.
If $x_3 x_4\in E(G)$, then
  by Lemma(b) $wx_3\in E(G)$ and thus,  by Lemma(b) as well, so is
   $wx\in E(G)$.
This is a contradiction. Thus $x_3 x_4\notin E(G)$. Then the edge
interchange which removes the edges $w x_4$ and $xx_3$ and inserts
the  edges $wx$ and $x_3 x_4$ produces a realization containing
the edge $wx$. By Case 1, $S$ has a realization containing a
$C_{k+1}$. This is a contradiction.\par
  {\bf Theorem\
2.} Let $m\ge 3$. Let $S$ be an n-term graphical sequence. Suppose
$S$ satisfies the following two conditions: (i) there is a
realization $G$ of $S$  containing a $C_{2m+1},$ such that for all
$ x, y\notin C_{2m+1}$, $d(x)=d(y)=m$ and $xy\notin E(G)$, (ii)
there is no realization of $S$  containing a $C_{2m+2}$. Then
$\sigma(S)\le m(2n-m-1)+2$.
\par {\bf Proof.} Let $C_{2m+1}$ be $w_1 w_2\cdots
w_{2m+1} w_1$, and $w_{2m+1+i}=w_i$. Since every realization of
$S$ does not contain a $C_{2m+2}$, by Lemma(a), for any $ v\notin
C_{2m+1}$, there is not $w_r, w_{r+1}$ such that $w_r v,
 w_{r+1} v\in E(G)$. Since for any $ x, y\notin C_{2m+1}$, $xy\notin E(G)$,
$d(x)=d(y)=m$, then $x, y$ are all adjacent to $m$ vertices of
$C_{2m+1}$.
 Assume without loss of generality  $w_1 x, w_4 x, w_6 x, \cdots, w_{2m}x\in E(G)$.\par
Case 1. \ Suppose there is $y\notin C_{2m+1}, y\ne x$ such that
there is a $w_i\in C_{2m+1}$ such that $w_i x\in E(G), w_i y\notin
E(G)$.
\par Subcase 1. \ Suppose $w_2 y\in E(G)$. By
Lemma(a), $w_3 y, w_1 y\notin E(G)$ and at most one vertex of
$w_4, w_5$ is adjacent to $y$.  If $w_6 y\in E(G)$, then $G$
contains a $C_{2m+2}:\ w_6 w_7\cdots w_{2m+1} w_1 x w_4 w_3 w_2 y
w_6$. This is a contradiction, thus $w_6 y\notin E(G)$. Next, if
$w_7 y\in E(G)$, then by Lemma(a), $w_8 y, w_6 y\notin E(G)$.
Since $y$ is adjacent to $m$ vertices of $C_{2m+1}$, Lemma(a)
forces $w_9 y, w_{11} y,\cdots,$ $ w_{2m+1} y\in E(G)$. Then $G$
contains a $C_{2m+2}:\ w_{2m+1} y$ $ w_2 w_1 x w_4 w_5\cdots
w_{2m} w_{2m+1}$. This is a contradiction, thus  $w_7 y\notin
E(G)$. Finally, suppose $w_6 y, w_7 y\notin E(G)$.
 Then, by Lemma(a),  $y$ at most is adjacent to
$m-1$ vertices of $C_{2m+1}$ - a contradiction.
\par Subcase 2. \ Suppose
$w_3 y\in E(G)$. By a similar method as Subcase 1 we can give a
contradiction.
\par Subcase 3. \ Suppose $w_2 y, w_3 y\notin E(G)$.
Lemma(a)  forces  $y$ to be adjacent to the following $m$ vertices
of $C_{2m+1}$: $w_1 , w_4 , w_6 , \cdots, w_{2m} $. This
contradicts the supposition of case 1.
\par Case 2. \ Suppose for any
$ y\notin C_{2m+1}, y\ne x,$ for any $ w_i\in C_{2m+1}$, if $w_i x
\in E(G)$, then $w_i y\in E(G)$. Then $w_1 y, w_4 y, w_6 y,
\cdots, w_{2m} y\in E(G)$.
\par Subcase 1. \ Suppose $w_2 w_5\in
E(G)$. Then $G$ contains $a$ $C_{2m+2}:\ w_5 w_2 w_3 w_4 x$ $ w_1
w_{2m+1} w_{2m}\cdots w_5$. This is a contradiction.
\par Subcase 2. \
Suppose $w_{2m+1} w_2\in E(G)$. Then $G$ contains a $C_{2m+2}:\
w_2 w_{2m+1} w_1 x w_{2m} w_{2m-1}\cdots w_2$. This is a
contradiction.
\par Subcase 3. \ Suppose there is an $i(3\le i \le m-1)$
such that $w_2 w_{2i+1}\in E(G)$. Then $G$ contains a $C_{2m+2}:\
w_{2i+1} w_2 w_3 w_4 \cdots w_{2i-2} xw_1 w_{2m+1} w_{2m}\cdots$ $
w_{2i+2} yw_{2i}$ $ w_{2i+1}$. This is a contradiction.
\par Subcase 4. \
Suppose $w_3 w_5\in E(G)$. Then $G$ contains a $C_{2m+2}:\ w_3 w_5
w_4 x w_6 $ $w_7 \cdots w_{2m+1} w_1 w_2 w_3$. This is a
contradiction.
\par Subcase 5. \ Suppose $w_3 w_{2m+1}\in E(G)$. Then
$G$ contains a $C_{2m+2}:\ w_{2m+1} w_3 w_2 w_1 x w_4 w_5 \cdots
w_{2m} w_{2m+1}$. This is a contradiction.
\par Subcase 6. \ Suppose
there is an $i(3\le i\le m-1)$ such that $w_3 w_{2i+1}\in E(G)$.
Then $G$ contains a $C_{2m+2}:\ w_{2i+1} w_3 w_2 w_1 x w_4
w_5\cdots w_{2i} y w_{2m} w_{2m-1}$ $\cdots w_{2i+1}$. This is a
contradiction.
\par Subcase 7. \ Suppose there is a $j (2\le j\le m-1)$
such that $w_{2j+1} w_{2m+1}\in E(G)$. Then $G$ contains a
$C_{2m+2}:\ w_{2j+1} w_{2m+1} w_{2m}\cdots w_{2j+2} x w_1 w_2
\cdots $ $w_{2j+1}$. This is a contradiction.
\par Subcase 8. \ Suppose
there is a $j$ and an $i (2\le j<i\le m-1)$ such that $w_{2j+1}
w_{2i+1}\in E(G)$. Then $G$ contains a $C_{2m+2}:\ w_{2j+1}
w_{2i+1} w_{2i}\cdots$ $w_{2j+2}$ $x w_{2i+2}$ $w_{2i+3}$ $\cdots$
$w_{2m+1} w_1 w_2\cdots w_{2j} w_{2j+1}$. This is a contradiction.
\par Subcase 9.  \ Suppose for any $i (\ 2\le i\le m), w_2
w_{2i+1}, w_3 w_{2i+1}\notin E(G)$, and for any $i, j ( 2\le j<
i\le m), w_{2j+1} w_{2i+1}\notin E(G)$.\par Then
$$\aligned &d(w_2),\ d(w_3)\le m+1\\
&d(w_5),\ d(w_7),\ \cdots,\ d(w_{2m+1})\le m\endaligned$$ Since
 for any  $y\notin C_{2m+1},\ d(y)=m$. Hence\par
$$\aligned \sigma(S)&\le m(n-2m-1)+d(w_2)+d(w_3)+d(w_5)+d(w_7)\\
&\hskip 12pt+\cdots+d(w_{2m+1})+d(w_1)+d(w_4)+d(w_6)+\cdots+d(w_{2m})\\
&\le m(n-2m-1)+2(m+1)+m(m-1)+(n-1)m\\
&=(n-2m-1+2+m-1+n-1)m+2\\
&=m(2n-m-1)+2.\endaligned$$ \par {\bf Theorem\ 3.} Let $m\ge 2$.
If $k=2m+1,\ n\ge 3m$, then $\sigma(C_k, n) =m(2n-m-1)+2$;\ if
$k=2m+2, n\ge 3m$, then $\sigma (C_k, n)\le m(2n-m-1)+2m+2$.\par
  {\bf Proof.} By [5]
theorem 2 and 3, $\sigma(C_5, n)=4n-4$ for $n\ge 5, \sigma (C_6,
n)=4n-2$ for $n\ge 7$. Clearly $\sigma(C_6, 6)=24$. Hence for
$m=2$, if $k=2m+1, n\ge 3m$, then $\sigma (C_k, n)\le
m(2n-m-1)+2$; if $k=2m+2, n\ge 3m$, then $\sigma (C_k, n)\le
m(2n-m-1)+2m+2$.
\par Suppose for t, $2\le t< m$, if $k=2t+1, n\ge
3t$, then $\sigma (C_k, n) \le t(2n-t-1)+2$ and if $k=2t+2, n\ge
3t$, then $\sigma (C_k, n)\le t(2n-t-1)+2t+2$.
\par Case 1. \ If $S$ is an n-term graphical sequence, with $k=2m+1, n\ge 3m
$,$\sigma(S)\ge m(2n-m-1)+2$.
  For $n=3m$,
$\sigma(S)\ge m(6m-m-1)+2=5m^2-m+2=2[{k-1 \choose 2} +{n-k+2
\choose 2} +1]$, which by [1] (chapter III, theorem 5.9) implies
that all realizations of $S$ contain a $C_k$. Now assume that
$S_1$ is a p-term graphical sequence, $3m\le p< n$,
$\sigma(S_1)\ge m(2p-m-1)+2$ and that there is a realization of
$S_1$ containing a $C_k$. We will show that if
$S=(d_1,d_2,...,d_{p+1})$ is a $p+1$-term graphical sequence with
realization $G$ and $\sigma(S)\ge m(2(p+1)-m-1)+2$, then  $S$ has
a realization containing a $C_{2m+1}$. Assume $d_1\ge
d_2\ge\cdots\ge d_n\ge 0.$
  Let $S'$ be the degree sequence of $G-v_{p+1}$
and suppose $d_{p+1}\le m$. Then $\sigma(S')\ge m(2(p+1)-m-1)+2-2m
=m(2p-m-1)+2$. Therefore, by our assumption, $S'$ has a
realization containing a $C_k$. Hence $S$ has a realization
containing a $C_k$. Thus, we may assume that $d_{p+1}\ge m+1$.
Since $\sigma(S)\ge m(2(p+1)-m-1)+2\ge
(m-1)(2(p+1)-(m-1)-1)+2(m-1)+2$, by our assumption, there is a
realization of $S$ containing a $C_{2m}$. Which by theorem 1
implies that $S$ has a realization containing a $C_{2m+1}$.
\par Case 2. \ If $k=2m+2$, $n\ge 3m,
S$ is an $n$-term graphical sequence with $\sigma(S)\ge
m(2n-m-1)+2m+2$ then we can prove, via a similar method as Case 1,
that $S$ has a realization containing a $C_{2m+2}$ .
\par Hence $\sigma(C_{2m+1}, n)\le
m(2n-m-1)+2$,  $\sigma(C_{2m+2}, n)\le m(2n-m-1)+2m+2$.\par By
[5], theorem 1 (Theorem A below), for $m\ge 2, k=2m+1, n\ge
2m+1,\sigma(C_k, n)\ge m(2n-m-1)+2$. Hence, for $m\ge 2$, if
$k=2m+1, n\ge 3m$, then $\sigma(C_k, n)=m(2n-m-1)+2. \hskip 30pt $
\par \vskip 10pt {\bf Lemma 4.} If $m\ge 3, n=3m+t(t=0,
1, 2, \cdots, 2m-2)$, then $\sigma(C_{2m+2}, n)$ $\le
m(2n-m-1)+2m+2-2[{t\over 2}]$.\par
 {\bf Proof.} By
theorem 3,  Lemma 4 holds for $t=0,1$. Now assume that Lemma 4
 holds for all $t-1, (1\le t\le 2m-2)$. We now prove that Lemma
holds for $t$.
\par  Let $S=(d_1,d_2,...,d_n)$ be an $n$-term graphical
sequence ($n=3m+t$), $G$  a realization of $S$ and $\sigma(S)\ge
m(2n-m-1)+2m+2 -2[{t\over 2}]$. Assume $d_1\ge d_2\ge \cdots \ge
d_n\ge 0$.\par   Let $S'$ be the degree sequence of $G-v_n$. If
$d_n\le m-1$, then $\sigma(S')\ge m(2n-m-1)+2m+2-2[{t\over 2}]
-2(m-1)\ge m(2(n-1)-m-1)+2m+2-2[{{t-1}\over 2}]$. By induction
suppose,
 $S'$ has a realization containing a $C_{2m+2}$. Hence $S$
has a realization containing a $C_{2m+2}$. Thus, we may assume
that $d_n\ge m$. Since $t\le 2m-2$, one has  $\sigma(S)\ge
m(2n-m-1)+2m+2-2[{t\over 2}] >m(2n-m-1) +2$. This implies, by
theorem 3, that $S$ has a realization containing a $C_{2m+1}$. Let
$w\in C_{2m+1}$, $x,y\notin C_{2m+1}$ and assume that every
realization of $S$  does not contain a $C_{2m+2}$. If  $d(x)\ge
m+1$, then since $d(w)\ge d_n\ge 3$, by theorem 1, $S$ has a
realization containing a $C_{2m+2}$. This is a contradiction.
Hence for any $ x\notin C_{2m+1}, d(x)=m$.
\par If
for any $ x, y\notin C_{2m+1}$, $xy\notin E(G)$, then, by theorem
2, $\sigma(S)\le m(2n-m-1)+2< m(2n-m-1)+2m+2-2[{t\over 2}]
\le\sigma(S)$. This is a contradiction. Thus, we may assume that
there is $x,y\notin C_{2m+1}$ such that $xy\in E(G)$. Let $S'$ be
degree sequence of $G-\{x, y\}$. Since $d(x)=d(y)=m$, then
$\sigma(S') \ge m(2n-m-1)+2m+2-2[{t\over
2}]-4m+2=m(2(n-2)-m-1)+2m+2-2[{{t-2}\over 2}]$. By induction
suppose, $S'$ has a realization containing a $C_{2m+2}$. Hence $S$
has a realization containing a $C_{2m+2}$. This is a
contradiction.\par Therefore $\sigma(C_{2m+2}, n)\le
m(2n-m-1)+2m+2-2[{t\over 2}].\hskip 30pt $\par  {\bf Theorem 5.}
 $\sigma(C_{2m+2}, n)=m(2n-m-1)+4$, for $m\ge 3,
n\ge 5m-2$.\par
  {\bf Proof.} By Lemma 4, for $m\ge
3, n=5m-2,\sigma (C_{2m+2}, n)\le m(2n-m-1)+2m+2-2[{{2m-2}\over
2}]=m(2n-m-1)+4$.\par Suppose for $p$, ( $5m-2\le p< n),\sigma
(C_{2m+2},\ p)\le m(2p-m-1)+4$. Let $S=(d_1,d_2,...,d_n)$ be an
$n$-term graphical sequence with realization $G$ and  $\sigma
(S)\ge m(2n-m-1)+4$. Assume $d_1\ge d_2\ge\cdots d_n\ge 0$ .
\par  If $d_n\le m$, then consider the degree
sequence,  $S^{\prime}$, formed by  $G-v_n$. Then $\sigma
(S^{\prime})$ $\ge m(2n-m-1)+4-2m=m(2(n-1)-m-1)+4$. By the
induction hypothesis, $S'$ has a realization containing a
$C_{2m+2}$. Hence $S$ has a realization containing a $C_{2m+2}$.
Thus, we may assume that  $d_n\ge m+1$. Since $\sigma (S)\ge
m(2n-m-1)+4\ge m(2n-m-1)+2$, theorem 3 implies that $S$ has a
realization containing a $C_{2m+1}$. Therefore, by theorem 1, $S$
has a realization containing a $C_{2m+2}$.\par Therefore $\sigma
(C_{2m+2}, n)\le m(2n-m-1)+4$.\par By [5] theorem 1 (Theorem A
below), for $m\ge 2, n\ge 2m+2, \sigma(C_{2m+2}, n)$ $\ge
m(2n-m-1)+4$. Hence $\sigma (C_{2m+2}, n)=m(2n-m-1)+4$ for $m\ge
3, n\ge 5m-2.\hskip 30pt$
\par
For completeness, we give a short proofs of the lower bounds for
$\sigma(C_{2m+1}, n)$ and $\sigma(C_{2m+2}, n)$ as following:
  \par
  {\bf Theorem A.} $\sigma(C_{2m+1}, n)\ge m(2n-m-1)+2$,\ for $n\ge
2m+1, m\ge 2$, $\sigma (C_{2m+2}, n) \ge m(2n-m-1)+4$, for $n\ge
2m+2, m\ge 2$.
\par
{\bf Proof.} By noting that  $G=K_m + \overline{K_{n-m}}$ gives a
uniquely realizable degree sequence and $G$ clearly does not
contain $C_{2m+1}$, $H=K_m + (\overline{K_{n-m-2}}\bigcup K_2)$
gives a uniquely realizable degree sequence and $H$ clearly does
not contain $C_{2m+2}$, this result can easily be seen. $\hskip
30pt$
 \section*{Acknowledgment}
 This paper was written in the University of Science and Technology of
China as a visiting scholar. The author thanks Prof. Li
Jiong-sheng for his advice. The author thanks the referees for
many helpful comments.
 \par

\end{document}